\documentclass{article}

\usepackage{amsmath}
\usepackage{amsfonts}
\usepackage{amscd}

\def\R{{\bf R}} 
\def\C{{\bf C}} 
\def\H{{\bf H}} 
\def\Z{{\bf Z}} 
\def\D{{\cal D}}

\newtheorem{theo}{Theorem}[section]
\newtheorem{lem}[theo]{Lemma}
\newtheorem{prop}[theo]{Proposition}
\newtheorem{cor}[theo]{Corolally}
\newtheorem{ex}[theo]{Example}
\newtheorem{rem}[theo]{Remark}
\newtheorem{defi}[theo]{Definition}

\begin{document}

\begin{center}
{\Large The first homology of the group of equivariant 
diffeomorphisms and its applications}
\end{center}

\begin{center} 
{\large
K\=ojun Abe
\footnote[2]{Department of Mathematical Sciences, Shinshu University, Matsumoto 390-6821, Japan

{\it E-mail address:} kojnabe@gipac.shinshu-u.ac.jp 

This research was partially supported by Grant-in-Aid for Scientific Research (No.16540058),
Japan Society for the Promotion of Science.}
and \ 
Kazuhiko Fukui 
\footnote[8]{Department of Mathematics, Kyoto Sangyo University, Kyoto 603-8555, Japan  \\ 
{\it E-mail address:} fukui@cc.kyoto-su.ac.jp 

This research was partially supported by Grant-in-Aid for Scientific Research (No.17540098),
Japan Society for the Promotion of Science.}}
\end{center}

%%%%%%%%%%%%%% Abstract %%%%%%%%%%%%%%
\bigskip
\noindent
{\bf\small Abstract}

\medskip
 Let $V$ be a representation space of a finite group $G$. 
 We determine the group structure of the first homology  
of the equivariant diffeomorphism group of $V$ . 
Then we can apply it to the calculation of the first homology 
of the corresponding automorphism groups of smooth orbifolds, compact Hausdorff foliations, 
codimension one or two compact foliations and a locally free $S^1$-action on the 3-sphere.
%%%%%%%%%%%%%%%%%%%%%%

\medskip



%%%%%%%%%%%%%%%%%%%% Section 1 %%%%%%%%%%%%%
%
\bigskip
\noindent
{\bf 1. \ Introduction}

 \setcounter{section}{1} 
  \setcounter{theo}{0}

\bigskip
 
  Let ${\cal D}(M)$ denote the group of diffeomorphisms of an $n$-dimensional smooth manifold $M$ which are 
 istopic to the identity through compactly supported diffeomorphisms.   
 Here smooth means differentiable of class $C^{\infty}$.
  In \cite{Th1}, Thurston proved that the group  ${\cal D}(M)$ is perfect, which means ${\cal D}(M)$ coincides 
  with its commutator subgroup. 
  It is known that the result is relevant to the foliation theory.   \par 
  
  When $M$ is a smooth manifold with boundary, Fukui \cite{F} proved that \\ 
  $H_1({\cal D}([0,1])) \cong \R \oplus \R$ and we proved in 
 \cite{A-F4} that the group ${\cal D}(M)$ is perfect if \\ 
 $n>1$.    
  Here the first  homology group of a group $K$ is given by $H_1(K)= K/[K,K]$.  We also have the same result in the relative versions.   
  There are many analogous results on the group of a smooth manifold $M$ 
  preserving a geometric structure of $M$.  \par 
   In this paper we first treat in the case of a representation space $V$ of 
  a finite group $G$. 
  Let ${\cal D}_G(V)$ denote the group of equivariant smooth diffeomorphisms 
of $V$ which are $G$-isotopic to the identity through compactly supported equivariant smooth diffeomorphisms. 
Then we shall prove that $H_1({\cal D}_G(V))$ is isomorphic to 
$H_1(Aut_G(V)_0)$ when $\dim V^G=0$ and ${\cal D}_G(V)$ is perfect when  
$\dim V^G>0$. 
Here $Aut_G(V)_0$ is the identity component of the group of $G$-equivariant 
linear automorphisms of $V$, and $V^G$ is the set of fixed points of 
$G$ on $V$.
It is easy to see that  $H_1(Aut_G(V)_0)$ is isomorphic to the group 
of the form ${\bf R}^d \times U(1)^{\ell}$. 
Then we can completely calculate $H_1({\cal D}_G(V))$ for any representation 
space $V$ of a finite group $G$. 
  
  Next we mention that the above result can be applied in various way. 
 First we apply it to the case of smooth orbifold. 
  Note that a smooth orbifold $N$ is locally diffeomorphic to the orbit space 
  $V/G$ of 
a representation space $V$ of a finite group $G$. Using the result by 
Biestone \cite{Bi1} and Schwarz \cite{Sc1}, we see that $H_1({\cal D}_G(V))$ 
is isomorphic to $H_1({\cal D}_G(V/G))$. 
 Combining the above result and the fragmentation lemma we can determine 
 the structure of 
 $H_1({\cal D}(N))$ of the diffeomorphism group ${\cal D}(N)$ for any smooth orbifold $N$. Then we see that 
 $H_1({\cal D}(N))$  describes a geometric 
 structure  around the isolated singularities.  \par 
 Let $M$ be a smooth $G$-manifold for a finite group $G$. Then 
 $H_1({\cal D}_G(M))$ is isomorphic to $H_1({\cal D}(M/G))$, and we see that  
 $H_1({\cal D}_G(M))$ describes the 
  properties  of the isotropy representations at the isolated fixed points 
 of $M$.  We remark that the above result can be also applied 
 to a smooth manifold $M$ with  properly disconnected smooth action 
 and also the modular group to detect the cusp points (see \cite{A}). \par 
   Secondly we apply it to the case of foliated manifold. Let ${\cal F}$ be 
 a compact Hausdorff foliation on a compact smooth manifold $M$. Let 
 ${\cal D}(M,{\cal F})$ denote 
 the group of all foliation preserving smooth diffeomorphisms of $(M,{\cal F})$ 
 which are isotopic to the identity through foliation preserving 
 diffeomorphisms. Then we can determine the structure of $H_1({\cal D} (M,{\cal F}))$ 
 using the result of Tsuboi (\cite {Ts}), and we see that it describes 
 the holonomy structures of isolated singular leaves of ${\cal F}$.  
   Calculating $H_1({\cal D} (M))$ for compact smooth orbifolds of 
 dimension less than 3 we shall determine the group structure of 
 $H_1({\cal D}(M,{\cal F}))$ 
 for compact foliated manifolds of codimension one and two. \par 
 
  Finally we apply it to a smooth G-manifold when $G$ is a 
 compact Lie group. If $M$ is a  principal $G$-manifold with $G$ 
 a compact Lie group, then we proved that the group ${\cal D}_G(M)$ 
 is perfect for  $\dim (M/G) > 0$ (Banyaga \cite {B1} and Abe and Fukui 
 \cite {A-F1}). In \cite {A-F2} we calculated $H_1({\cal D}_G(M))$ 
 when $M$ is a smooth $G$-manifold with codimension one orbit. 
 We shall apply the above result to the case of a locally free 
 $U(1)$-action on  the 3-sphere, and calculate $H_1({\cal D}_{U(1)}(S^3))$. 

 The paper is organized as follows. In $\S 2$, we state 
 the main theorem and apply it to the definite cases. 
 $\S 3$ is devoted to prove the main theorem.  In $\S 4$, we apply 
 the main result to the group of diffeomorphisms of an orbifold and 
 calculate the first homology of the group. In \S 5 and \S 6, we apply 
 the main result to compact Hausdorff foliations and calculate 
 the first homology of the groups of foliation preserving diffeomorphisms. 
 In $\S 7$, we apply the result to the case of a locally free 
 $S^1$-action on the 3-sphere and calculate the first homology of 
 the group of equivariant diffeomorphisms of $S^3$. 
%%%%%%%%%%%

The authors would like to thank Takashi Tsuboi who kindly answered their query 
on Theorem 3.1.

%
%%%%%%%%%%%%% Section 2 %%%%%%%%%
%
\bigskip
\newpage
\noindent
{\bf 2. \ The first homology of ${\cal D}_G(V)$ and applications 
to the low dimensional cases.}

 \setcounter{section}{2} 
  \setcounter{theo}{0}

\bigskip 

 Let $V$ be an $n$-dimensional real 
representation space of a finite group $G$.
Let ${\cal D}_G(V)$ be the subgroup of ${\cal D}(V)$ 
which consists of $G$-equivariant smooth diffeomorphisms of $V$ 
which are isotopic to 
the identity through $G$-equivariant smooth isotopies 
with compact support. 

In this section we shall state the main theorem with respect to the first homology $H_1({\cal D}_G(V))$ and calculate some definite cases. 
Let $Aut_G(V)$ be the group of $G$-equivariant automorphisms of $V$ and 
$Aut_G(V)_0$ its identity component. Let $V^G$ be the set of fixed points of 
$G$ on $V$. 
Let ${\mathfrak D}:\ \D_G(V) \to A_G(V)_0$ be the map defined by 
${\mathfrak D}(f) = df(0).$

Then we have the following, which is our main result. 

%
%%%%%%%%%%% Th 2.1 %%%%%%
%
\begin{theo}\label{th2.1} \quad 
$(1)$ If $\dim V^G > 0$, then ${\cal D}_G(V)$ is perfect.

\noindent
$(2)$ If $\dim V^G = 0$, then 
${\mathfrak D}_*: H_1(\D_G(V)) \to H_1(A_G(V)_0)$ is isomorphic. 
\end{theo}
%
%%%%%%%%%
%

The representation $V$ of $G$ is expressed as 
 $\displaystyle V=\oplus_{i=1}^d k_iV_i$, where $V_i$ runs over the inequivalent irreducible representation spaces of $G$ and $k_i$ is a positive integer. Let $Hom_G(V_i)$ be the set of $G$-equivariant homomorphisms of $V_i$. 
%
%%%%%%%%%%%%%%% Cor 2.2 %%%%%%%%%%%%
%
\begin{cor}\label{cor2.2} \quad 
If $\dim V^G = 0$, then 
 $H_1({\cal D}_G(V)) \cong {\bf R}^d \times \overbrace{U(1) \times 
\cdots \times U(1)}^{d_2}$, 
where $d_2$ is the number of $V_i$ with $\dim Hom_G(V_i)= 2 $.
\end{cor}
%
%%%%%%%%%%%%
%

\noindent
{\bf Proof}.\quad 
 If $V_i$ is a real restriction of an irreducible complex 
 (resp. quaternionic) representation of $G$, then $Hom_G(V_i)$ is 
 isomorphic to $\C$ (resp. $\H$). Otherwise $Hom_G(V_i)$ is isomorphic to 
 $\R$ (c.f. Adams \cite{Ad}, 3.57).  
Then  
$$Aut_G(k_i V_i) \cong \left\{\begin{array}{lll}
GL(k_i, \R) & \mbox{if} \quad   Hom_G(V_i) \cong \R \\
GL(k_i, \C) & \mbox{if} \quad  Hom_G(V_i) \cong \C \\
GL(k_i, \H) & \mbox{if} \quad  Hom_G(V_i) \cong \H. 
\end{array}\right.$$
Note that $GL(k_i, F) \cong SL(k_i, F) \times F$ for $F = \R,\, \C,\, \H$, 
where $SL(k_i, F)$ is the special linear group. Since  $SL(k_i, F)$ is 
a simple group,  
$$H_1(Aut_G(k_i V_i)) \cong \left\{\begin{array}{ll}
\C^*   & \mbox{if} \quad   Hom_G(V_i) \cong \C \\
\R^* & \mbox{if} \quad   Hom_G(V_i) \cong \R \ \mbox{or} \ \H. \\ 
\end{array}\right.$$
Then 
$$H_1(Aut_G(k_i V_i)_0) \cong \left\{\begin{array}{ll}
U(1) \times \R   & \mbox{if} \quad  \dim Hom_G(V_i) =2 \\
\R & \mbox{if} \quad  \dim Hom_G(V_i) = 1 \ \mbox{or} \ 4. \\ 
\end{array}\right.$$
Combining Theorem \ref{th2.1} with Schur's lemma we have 
$$H_1(Aut_G(V)_0) \cong H_1(Aut_G(k_1 V_1)_0) \times \cdots \times 
 H_1(Aut_G(k_d V_d)_0)$$
 $$\cong {\bf R}^d \times \overbrace{U(1) \times 
\cdots \times U(1)}^{d_2}. $$
This completes the prooof of Corollay \ref{cor2.2}.

 \bigskip
Now we describe the following Examples. 

%
%%%%%%%%%%%%% Ex 2.3 %%%%%%
%
\begin{ex} \label{ex2.3} \quad
 Let $V$ be the $n$-dimensional representation space of the cyclic group 
${\bf Z}_2$ with $\dim V^G = 0$. Then 
$H_1({\cal D}_{{\bf Z}_2}(V)) \cong {\bf R}$.
\end{ex}
%
%%%%%%%%%%%%%
%

{\bf Proof}.\quad 
Let $\tilde{\R}$ be the non-trivial one dimensional representation space 
of $\Z_2$. 
Then $V \cong n \tilde{\R}$. The proof follows from Theorem \ref{th2.1} (2) 
and Corollary \ref{cor2.2}. 

\bigskip
Next, we consider the case $V$ is a non-trivial 2-dimensional 
orthogonal representation space 
of a finite subgroup of $G$. Then the finite subgroups $G$ of $O(2)$ are 
classified as follows. 
$G$ is isomorphic to either ${\bf Z}_k \ (k \ge 2)$ which acts on 
$V$ as $k$ rotations, the reflection group ${\bf D}_1$ which acts on $V$ as a reflection or the dihedral group ${\bf D}_l=\{ u,v \,; u^l=v^2=(uv)^2=1 \} \ (l \ge 2)$
which acts on $V$ as $l$ rotations and $l$ reflections.  
Then we have the following. 
  
%
%%%%%%%%%%%%% Ex 2.4 %%%%%%%
%
\begin{ex} \label{ex2.4} \quad
 $H_1({\cal D}_G(V)) \cong \left\{\begin{array}{llll}
{\bf R} & \mbox{if}\ G \cong {\bf Z}_2 \\
{\bf R}\times U(1) & \mbox{if}\ G \cong {\bf Z}_n \ (n \ge 3) \\
0 & \mbox{if}\ G \cong {\bf D}_1 \\
{\bf R}^2 & \mbox{if}\ G \cong {\bf D}_2  \\
{\bf R} & \mbox{if}\ G \cong {\bf D}_n \ (n \ge 3). 
\end{array}\right.$ 
\end{ex}
%
%%%%%%%%%%%%%
%

{\bf Proof}.\quad 
If $G \cong {\bf D}_1$, then $\dim V^G =1$.  
It follows from Theorem 2.1 (1) that the group ${\cal D}_G(V)$ is perfect.  
 For the case $G \cong {\bf Z}_n \, (n \ge 3)$, $V$ is the real restriction of 
 one dimensional irreducible complex representation of ${\bf Z}_n$. 
 By Theorem 2.1 (1) and Corollary 2.2, $H_1({\cal D}_G(V)) \cong {\bf R}\times U(1)$. 
If $G \cong {\bf Z}_2$, then it is a special case of Example 2.3. 
 If $G \cong {\bf D}_n \, (n \ge 3)$, then 
$V$ is a 2-dimensional real irreducible representation of $G$ 
which is not isomorphic to the real restriction of a complex representation 
of $G$.  Then it follows from Theorem 2.1 (2) and Corollary 2.2 that 
$H_1({\cal D}_G(V)) \cong {\bf R}$. 
If $G \cong {\bf D}_2$, then $V$ is expressed as the direct sum of two 
 inequivalent 1-dimensional real representations 
of $G$. Then it follows from Theorem 2.1 (2) and Corollary 2.2 that 
$H_1({\cal D}_G(V)) \cong {\bf R}^2$.
 This completes the proof. 

\vspace{5mm} 

%
%%%%%%%%%%%%%%%%% Section 3 %%%%%%%
%
%\newpage

\noindent
{\bf 3. \ Proof of Theorem 2.1.}

\setcounter{section}{3} 
  \setcounter{theo}{0}

\bigskip

In this section we shall prove Theorem \ref{th2.1}. 
For this purpose we describe the result on the group of leaf preserving 
diffeomorphisms of the product foliation. \par 
 If $V$ is an $n$-dimensional representation space of $G$ with $\dim V^G =n-p$, 
 then $V$ can be expressed as 
 $V=W \times {\bf R}^{n-p}$ such that $W$ is a representation space of $G$ 
with $\dim W^G=0$.  
 Let ${\cal F}_0$ be the product foliation of $V$ with leaves of the 
 form $\{ \{ x \} \times {\bf R}^{n-p} \}$, where $(x,y)$ is a coordinate of 
 $V=W \times {\bf R}^{n-p}$. 
By ${\cal D}_L(V, {\cal F}_0)$ we denote the group of leaf preserving 
smooth diffeomorphisms of $(V, {\cal F}_0)$ which are isotopic to the 
identity through leaf preserving smooth diffeomorphisms with compact support. 
Let ${\cal D}_{L,G}(W \times {\bf R}^{n-p}, {\cal F}_0)$ denote 
the subgroup of ${\cal D}_L(V, {\cal F}_0)$  consisting of leaf preserving 
$G$-equivariant smooth diffeomorphisms of 
$(W \times {\bf R}^{n-p}, {\cal F}_0)$ which are isotopic to the identity 
through leaf preserving $G$-equivariant smooth diffeomorphisms 
with compact support. 

T.Tsuboi \cite{Ts} proved the perfectness of ${\cal D}_L(V, {\cal F}_0)$ by looking at the proofs in \cite {H} and \cite {Th1}. Futhermore he observed 
the following equivariant version. 
The point of the proof is observing the group  
${\cal D}_{L,G}(W \times T^{n-p} , {\cal F}'_0)$ to be perfect and 
using the equivariant fragmentation lemma. Here ${\cal F}'_0$ denotes 
the product foliation of $W \times T^{n-p}$ with leaves of the form 
$\{ pt \} \times T^{n-p}$.

%
%%%%%%%%%%% Th 3.1 %%%%%%
%
\begin{theo}\label{th3.1} {\rm (T.Tsuboi \cite{Ts})} \quad 
${\cal D}_{L,G}(W \times {\bf R}^{n-p}, {\cal F}_0)$ is perfect.
\end{theo}
%
%%%%%%%%%
%

\bigskip

We shall prove Theorem \ref{th2.1} (1) by induction of the order $|G|$ 
of $G$. Note that, from the 
result by Thurston \cite{Th1}, it holds when $|G| = 1$. 

%
%%%%%%%%%%% Th 3.2 %%%%%%
%
\begin{theo}\label{th3.2} 
${\cal D}_G (W^p \times {\bf R}^{n-p})$ is perfect when $n-p \ge 1$.
\end{theo}
%
%%%%%%%%%
%

\noindent
{\bf Proof.}\quad
We assume that Theorem \ref{th3.2} holds for any finite subgroup $H$ 
with $|H| < |G|$.  
For each $f \in {\cal D}_G (W \times {\bf R}^{n-p})$ we shall prove that 
$f$ can be written as the composition of commutators. 
The proof is divided into two steps. \par 
{\it First step.} \quad We shall prove that $f$ can be expressed as 
$f= f_1 \circ f_2$ 
such that \par 
 (a)\ $f_1 \in [{\cal D}_G (W \times {\bf R}^{n-p}),\ 
{\cal D}_G (W \times {\bf R}^{n-p})]$ \par 
 (b)\  $f_2 \in {\cal D}_G (W \times {\bf R}^{n-p})$ with $supp(f_2) \subset 
W \times {\bf R}^{n-p} \setminus \{0\} \times {\bf R}^{n-p}$.  \\
For the purpose, by using the fragmentation lemma (c.f. \cite{A-F1}, Lemma 1), 
we can assume that $f$ is close to the identity in the $C^{\infty}$-topology 
 of ${\cal D}_G (W \times {\bf R}^{n-p})$ and 
 $supp(f)$ is contained 
 in the $\delta$-neighborhood  $B_{\delta }$ of a point 
 $(0,0) \in \{0\} \times {\bf R}^{n-p}$ for some $0 <\delta < 1$.
  Then, there exist $g_1, g_2 \in {\cal D}_G (W \times {\bf R}^{n-p})$ 
  such that 

\vspace{1mm}
$(1)$ $f=g_2 \circ g_1$, \par 
$(2)$ $g_1$ and $g_2$ are close to the identity supported in $B_{\delta }$, 
\par 
 $(3)$ $g_1 \in {\cal D}_{L,G}(W \times {\bf R}^{n-p}, {\cal F}_0)$, \par 
 $(4)$ 
 $g_2$ is written of the form $g_2(x,y) = (\hat{g}_2(y)(x),y)$ such that 
 $\hat{g}_2(y) \in {\cal D}_{G}(W)$ for any 
 $x \in W,\ y \in {\bf R}^{n-p}$.

 \vspace{1mm}
\noindent
By Theorem \ref{th3.1}, $g_1$ is expressed as a product of commutators of elements in ${\cal D}_{L,G}(W \times {\bf R}^{n-p}, {\cal F}_0)$ which is the 
 subgroup of ${\cal D}_G (W \times {\bf R}^{n-p})$. 

%%%%%%%%%%%%%
Next we shall prove that $g_2$ can be expressed as a product of commutators of 
elements in ${\cal D}_G (W \times {\bf R}^{n-p})$. 
Let d$\hat{g_2}(y)(0) \in Aut_G(W)_0$ be the differential of 
$\hat{g_2}(y)$ at $x=0$ for each $y \in {\bf R}^{n-p}$. 
Since $g_2$ is close to  the identity in the $C^{\infty}$-topology, 
the map $d \hat{g_2}(\cdot)(0):\ \R^{n-p} \to Aut_G(W)_0$ is $C^1$-close to 
the unit map $1_{W}$ and is supported in $B_{\delta }'$, where 
$B_{\delta }'$ is the $\delta$-neighborhood at 0 in $\R^{n-p}$.
Let $q=\dim Aut_G(W)_0$.  Then it follows from Lemma 4 of \cite {A-F1} 
that there exist 
$G_i : {\bf R}^{n-p} \to Aut_G(W)_0$ and 
$\varphi _i \in {\cal D}({\bf R}^{n-p})$ \, 
$(i=1,2,\cdots ,q)$,  satisfying that 

\vspace{1mm}
(i)  $\varphi _i$ is $C^1$-close to the identity and is suppoted in 
$B_{\delta }'$,  \par 
(ii) $G_i$ has a compact support and the image of $G_i$ is contained 
in a sufficiently small neighborhood of $1_W$ and \par 
(iii)  
$d \hat{g_2}(\cdot)(0) =(G_1^{-1}\cdot(G_1\circ \varphi _1))\cdot \cdots \cdot(G_q^{-1} 
\cdot(G_q\circ \varphi _q)).$

\vspace{1mm}

Let $\mu : W \to [0,1]$ be a $G$-equivariant smooth 
function satisfying that $\mu(x)=1$ for $\displaystyle ||x|| \le \frac12 $ 
and $\mu(x)=0$ for $||x|| \ge 1$ and $\mu$ is monotone decreasing with 
respect to the distance from the origin. 
Then we have smooth diffeomorphisms  
$h_{G_i }: W \times {\bf R}^{n-p} \to W \times {\bf R}^{n-p} \ (i=1,...,q)$ 
defined by 
$$h_{G_i}(x,y)=\left( \mu(x) G_i(y)(x)+(1-\mu(x)) x, \,y \right).$$ 
By the property (ii), $h_{G_i}$ is a $G$-equivariant smooth diffeomorphism of 
$W \times {\bf R}^{n-p}$ and $h_{G_i}(x,y)=(G_i(y)x,y)$ for 
$\displaystyle ||x|| \le \frac12 $. 
For each $\varphi_i \in {\cal D}({\bf R}^{n-p})$, we put 
$$F_{\varphi_i}(x,y)=(x, \mu(x)\varphi_i(y)+(1-\mu(x))y).$$ 

\noindent
Since $\varphi_i$ is $C^1$-close to the identity, 
 $F_{\varphi_i}$ is a $G$-equivariant smooth diffeomorphism of 
 $W \times {\bf R}^{n-p}$ . 
 If $\displaystyle ||x|| \le \frac12$, then 
  $F_{\varphi_i}^{-1}(x,y)=(x,\varphi^{-1}_i(y))$. Then we have that 
$$h_{G_i }^{-1}\circ F_{\varphi_i}^{-1} \circ h_{G_i }\circ F_{\varphi_i}(x,y)
=( G_i (y)^{-1}\cdot G_i (\varphi_i(y))\cdot x, y).$$ 
for $\displaystyle ||x|| \le \frac12$.

 Now we obtain the following by easy calculations. 
%
%%%%%%%%% Lemma 3.3 %%%%%%
%
\begin{lem}\label{lem3.3} \quad  On a neighborhood of 
\,$0 \times {\bf R}^{n-p}$,  
$\displaystyle h_{{\rm d}\hat{g_2}(\cdot)(0)}$ is coincide with \\
$ \prod_{i=1}^q[h_{G_i }^{-1}, F_{\varphi_i}^{-1}]$ which is contained in 
$[{\cal D}_G (W \times {\bf R}^{n-p}), {\cal D}_G (W \times {\bf R}^{n-p})]$. 
\end{lem}
%
%%%%%%%%%%
%

Put $g_3=(\prod_{i=1}^q[h_{G_i }^{-1}, F_{\varphi_i}^{-1}])^{-1}\circ g_2$. 
Then $g_3$ has of the form $g_3=(\hat{g}_3(y)(x), y)$ 
such that 
 $\hat{g}_3(y) \in {\cal D}_{G}(W)$ for any  $x \in W,\ y \in {\bf R}^{n-p}$. 
 Note that $g_3$ is a $G$-equivariant 
smooth diffeomorphism and the differential $d\hat{g}_3(y)(0)= 1_W$  
for any $y \in \R^{n-p}$ and $supp(\hat{g}_3) \subset B'_{\delta}$. \par 
Let $c$ be a real number with $0<c<1$. 
 Let $\psi \in {\cal D}_G (W \times {\bf R}^{n-p})$ such that   
 $\psi (x,y)=(cx,y)$ for $y \in B'_{\delta}$ 
and $supp(\psi) \subset B_{2\delta}$. 

 Then the Jacobi matrix $J((g_3\circ \psi) (\cdot\, ,\, y))$ at $x=0$ is 
 the scalar matrix by $c$  for each $y \in B'_{\delta}$. 
Note that each diffeomorphism $(g_3\circ \psi) (\cdot\, ,\, y)$ satisfies 
the conditions in Sternberg \cite{S1} and \cite {S2}. By using the parameter 
version of Borel theorem (c.f. Narashimham \cite{n}, \S 1.5.2), 
 we have the parameter version of Sternberg's theorem.  
 Then there exists a smooth diffeomorphism $R$ of $W \times {\bf R}^{n-p}$ 
 with compact support such that \par 
 (1)\  $(R^{-1}\circ (g_3\circ \psi)\circ R)(x,y) = \psi(x,y)$,  
 for $(x,y) \in B_{\delta}.$  \par 
 (2)\  $R(x,y)=(x+ S(x,y),y)$,  \\
 where $S$ is $1$-flat 
 at each point of $\{ 0 \} \times {\bf R}^{n-p}$, as a mapping of $x$, with compact support. \par
   By averaging on the group $G$, we define  
 $$\hspace{3mm} \tilde{R}_t(x,y)
  =\displaystyle \frac{1}{|G|} \sum_{g \in G}(x+ g^{-1}\cdot 
 (t S(g\cdot x,y)),\ y) $$  
for $(x,y) \in W \times {\bf R}^{n-p},\ 0 \le t \le 1$. 
Here $|G|$ is the order of $G$. 
Then  $\tilde{R}_t :  W \times {\bf R}^{n-p} \to  W \times {\bf R}^{n-p}$ 
is a smooth $G$-map. 
Since $S$ is $1$-flat  at each point of $\{ 0 \} \times {\bf R}^{n-p}$, 
 as a mapping of $x$, $\{\tilde{R}_t\}$ is a smooth $G$-isotopy on 
a neighborhood of $\{ 0 \} \times {\bf R}^{n-p}$ 
in $W \times {\bf R}^{n-p}$.  
Note that $\{\tilde{R}_t\}$ has a compact support and $\tilde{R}_0$ is 
the identity map. 
Combining the isotopy extension theorem with the isotopy integration 
theorem by Bredon \cite{Br}, Chapter IV, Theorem 3.1, 
there exists a $G$-isotopy 
$\{\hat{R}_t\}$ with compact support such that $\{\hat{R}_t\}$ 
coincides with $\{\tilde{R}_t\}$ on a neighborhood $U$ 
of $\{ 0 \} \times {\bf R}^{n-p}$ and  $\hat{R}_0$ is the identity map.  \par  
If $(x,y) \in U$, then  
\begin{eqnarray*}
\psi\circ \hat{R}_1(x,y) & = & \displaystyle \psi 
(\frac{1}{|G|} \sum_{g \in G} g^{-1}\cdot R(g\cdot (x,y))) \\
& = & \displaystyle \frac{1}{|G|} \sum_{g \in G} g^{-1}\cdot 
\psi(R(g\cdot (x,y))) \\
& = & \displaystyle \frac{1}{|G|} \sum_{g \in G} g^{-1}\cdot 
(R\circ g_3 \circ \psi(g\cdot (x,y)) \\
& = & \displaystyle \frac{1}{|G|} \sum_{g \in G} g^{-1}\cdot 
R(g\cdot (g_3 \circ \psi(x,y))) \\
& = & \hat{R}_1(g_3 \circ \psi(x,y)).
\end{eqnarray*}
Thus $\psi\circ \hat{R}_1=\hat{R}_1\circ g_3 \circ \psi$ on $U$.  
Put 
$$g_4=g_3 \circ (\hat{R}_1^{-1}\circ \psi\circ \hat{R}_1\circ 
\psi^{-1})^{-1}.$$
Then $g_4=1$ on $U$ and $g_3^{-1} \circ g_4$ is containted in 
the commutator subgroup of ${\cal D}_G (W \times {\bf R}^{n-p})$. \par 
  Therefore we have proved that each 
  $f \in {\cal D}_G (W \times {\bf R}^{n-p})$ 
can be written as $f= f_1 \circ f_2$ satisfying the conditions (a),(b). \par 
 {\it Second step.} \quad 
   We shall prove $f_2$ is written as the composition of commutators 
 in ${\cal D}_G (W \times {\bf R}^{n-p})$. 
 There exist finite points $ p_i \in W \times {\bf R}^{n-p} $ 
  $(i=1,\cdots , k)$ and open disk neighborhoods $U(p_i)$ at $p_i$ 
  $(i=1,\cdots , k)$ such that each $U(p_i)$ is a linear slice at $p_i$ and 
  $supp(f_2) \subset \displaystyle \bigcup_{i=1}^k G\cdot U(p_i)$. 
By the fragmentation lemma (see  \cite {A-F1}, Lemma 1), 
there are $h_i \in {\cal D}_G(W \times {\bf R}^{n-p})$ $(i=1,\cdots , \ell)$ 
such that 

(1) each $h_j$ is isotopic to the identity through $G$-diffeomorphisms 
with support in $G\cdot U(p_i)$, 

(2) $f_2=h_1\circ \cdots \circ h_{\ell}$. 

\noindent
Since $U(p_i)$ is a linear slice at $p_i$, the isotropy subgroup $G_{p_i}$ 
acts linearly 
on $U(p_i)$ and $G\cdot U(p_i)$ is a disjoint union of $|G/G_{p_i}|$ disks. 
Then from the condition (1), we have 
$$h_i(g\cdot U(p_i))=g\cdot U(p_i)\ \, \mbox{for}\ \, g\in G.$$ 
Thus $h_i$ is determined by the restriction map 
$h_i|_{U(p_i)} \in {\cal D}_{G_{p_i}}(U(p_i))$. 

 From the assumption of induction, each $h_i|_{U(p_i)}$ can be written 
 as a product of 
commutators of elements in ${\cal D}_{G_{p_i}}(U(p_i))$. 
Then  each $h_i$ is contained in the commutator subgroup of 
 ${\cal D}_G(W \times {\bf R}^{n-p})$. 
This completes the proof of Theorem \ref{th3.2}.  

\bigskip

\noindent
{\bf Proof of Theorem 2.1 (2).} \
 Let $V$ be an $n$ dimensional representation space of a finite group $G$ 
 with $\dim V^G=0$. 
Let $\Phi : {\cal D}_G(V) \to Aut_G(V)_0$ be the homomorphism defined by 
 $\Phi(f)=df(0)$ for any $f \in {\cal D}_G(V)$. 
 It is easy to see that $\Phi$ is epimorphic. 

Since 
$\begin{CD} 1 \to \ker \Phi \stackrel{\iota}{\to} {\cal D}_G (V) 
\to Aut_G(V)_0 \to 1 \end{CD}$ is exact, 
we have the following exact sequence: 
$$\begin{CD}
\ker \Phi/[\ker \Phi, {\cal D}_G (V)] \stackrel{\iota_*}{\to} 
H_1({\cal D}_G (V)) \to H_1( Aut_G(V)_0) \to 1.
\end{CD}$$
  The proof of Theorem $\ref{th2.1}\, (2)$ completes from the following. 
%
%%%%%%%%% Proposition 3.4 %%%%%%
%
\begin{prop}\label{prop3.4} \quad   
$\ker \Phi = [\ker \Phi,\ {\cal D}_G (V)]$ .
\end{prop}
%
%%%%%%%%%%
%
  Let $f \in \ker \Phi$. 
  We shall prove that $f \in [\ker \Phi,\ {\cal D}_G (V)]$. 
    We can assume that 
  $supp(f) \subset B_{\delta}$, where $B_{\delta}$ is the $\delta$-neighborhood of $0$ in $V$.  
 Let $c$ be a real number with $0<c<1$. 
 Let $\psi \in {\cal D}_G (V)$ such that   
 $\psi (x)= cx$ for $x \in B_{\delta}$ 
and $supp(\psi) \subset B_{2\delta}$. 
Then the Jacobi matrix $J(f \circ \psi)$ at $x=0$ is  the scalar matrix 
by $c$. 
By the parallel argument as in the proof of Theorem \ref{th3.2}, we can find 
$\hat{R}_1 \in {\cal D}_G (V)$ such that $\hat{R}_1$ is 1-tangent to 
the identity at the origin 
and $\psi\circ \hat{R}_1=\hat{R}_1\circ f \circ \psi$ 
on a neighborhood $U$ of $0$ in $V$. 
Put 
$f_2 = f \circ (\hat{R}_1^{-1}\circ \psi\circ \hat{R}_1\circ \psi^{-1})^{-1}.$
Then $f_2=1$ on $U$. \par 
   There exist finite points $ p_i \in V-U $  
  $(i=1,\cdots , k)$ and open disk neighborhoods $U(p_i)$ at $p_i$ 
  $(i=1,\cdots , k)$ such that each $U(p_i)$ is a linear slice at $p_i$ and 
  $supp(f_2) \subset \displaystyle \bigcup_{i=1}^k G\cdot U(p_i)$ 
  such that 

(1) each $h_i$ is isotopic to the identity through $G$-diffeomorphisms \\ 
with support in $G\cdot U(p_i)$, 

(2) $f_2=h_1\circ \cdots \circ h_{\ell}$.  \par 
  Since $\dim U(p_i)^{G_{p_i}} >0$, by Theorem \ref{th3.2} 
each $h_i| U(p_i)$ is contained in the commutator subgroup of 
$D_{G_{p_i}} (U(p_i))$.  
Then  each $h_i$ is contained in the commutator subgroup  
 $[\ker \Phi,\ {\cal D}_G (V)]$. 
This completes the proof of Proposition \ref{prop3.4}.

\vspace{5mm}

%%%%%%%%
\bigskip
%\newpage
\noindent
{\bf 4. Application to orbifolds.}

 \setcounter{section}{4} 
  \setcounter{theo}{0}

\bigskip
In this section we study the groups of smooth diffeomorphisms of 
smooth orbifolds. 

\bigskip

%
%%%%%%%%%%%%%% Def 4.1 %%%%%%%%%%%
%
\begin{defi}\label{def4.1}
{\rm (\cite {Sa} and \cite {Th2})}.\, \quad 
A paracompact Hausdorff space $M$ 
is called a smooth orbifold if there exists an open covering $\{ U_i \}$ 
of $M$, 
closed under finite intersections, satisfying the following conditions. 
\par
$(1)$ For each $ U_i $, there are a finite group $\Gamma_i$, 
a smooth effective action of $\Gamma_i$ 
on an open set $\tilde{U_i} $ of ${\bf R}^n$ and a homeomorphism 
$\phi_i : U_i \to \tilde{U_i}/\Gamma_i$. We call $\tilde{U_i} $ a cover of 
$U_i$.

$(2)$ Whenever $U_i \subset U_j$, there is a smooth imbedding 
$\phi_{ij} : \tilde{U_i} \to \tilde{U_j}$ such that the following 
diagram commutes:

$$\begin{CD}
\tilde{U_i} @>\phi_{ij}>>\tilde{U_j} \\
@V\phi_i^{-1}\circ\pi_iVV @VV\phi_j^{-1}\circ\pi_jV \\
U_i @>>> U_j, 
\end{CD}$$
where $\pi_k : \tilde{U_k} \to \tilde{U_k}/\Gamma_k (k=i,j)$ are 
the natural projections. 
Each $(U_i,\phi_i)$ is called a local chart of $M$.
\end{defi}
%
%%%%%%%%%%%%%%%%%%%%%%%%%%%%%%%%%%%%%%%%%%%%%%%%%%%%%%%%%%%%%%%
%

\vspace{2mm}

 Now we consider smooth diffeomorphisms on a smooth orbifold. 
 We refer to Bierstone \cite {Bi1}, \cite{Bi2} and Schwarz \cite{Sc1}, 
 \cite{Sc2}. 
 Let $M$ be a smooth orbifold. A continuous function $h : M \to \R$ is 
said to be smooth if for any local chart 
$(U_i, \phi_i)$ of $M$ 
 the composition $h \circ \phi_i^{-1} \circ \pi_i$ is smooth. 
A continuous map $f: M \to M$ is said to be smooth if  for any smooth function 
$h: M \to \R$, the composition $h \circ f$ is smooth. A homeomorphism 
 $f : M \to M$ is called a smooth diffeomorphism if 
$f$ and $f^{-1}$ are smooth maps. 
Let ${\cal D} (M)$ denote the group of smooth diffeomorphisms of $M$ 
which are isotopic to the identity through smooth diffeomorphisms 
with compact support.  \par 

 Let $(U_i,\phi_i)$ and $(U_j,\phi_j)$ be local charts of $M$ such that 
 $U_i$ and $U_j$ are diffeomorphic. We can assume that $\tilde{U_i}$ 
 and $\tilde{U_j}$ are invariant open neighborhoods of the origin of 
 representation spaces of $\Gamma_i$ and $\Gamma_j$ respectively.  
 Then by the result of Strub \cite{St}, $\Gamma_i$  and $\Gamma_j$ 
 are isomorphic 
 and the corresponding representations are equivalent. 
  Let $(\phi_i^{-1} \circ \pi_i)_*:\  \tilde{\cal D}_{\Gamma_i} (\tilde{U}_i) 
 \to {\cal D} (U_i)$ be the natural homomorphism. 
 Since $(\phi_i^{-1} \circ \pi_i)_*$ is epimorpic by  Bierstone \cite{Bi1}, 
 Theorem B, we can give the induced tpology on ${\cal D} (U_i)$ . 
 We can naturally regard ${\cal D} (U_i)$ as a 
subgroup of ${\cal D} (M)$.  
Then we give the topology on $\D(M)$ such that a subset $O$ of $\D(M)$ is 
open if $ O \cap \D(U_i)$ is open in $\D(U_i)$ for each local chart 
$(U_i,\phi_i)$. 
   
   \bigskip

%
%%%%%%%%%%%%%% Def 4.2 %%%%%%%%%%%%
%
\begin{defi}
\quad A point $x \in M$ is said to be an 
{\it isolated singular point} 
if for a local chart $(U_i,\phi_i)$  around $x$, the group $\Gamma_i$ 
acts on $\tilde{U_i}$ such that $\tilde{x}$ is the isolated fixed point, 
where $\tilde{x} \in \tilde{U_i}$ with $\pi_i(\tilde{x})=x$.
By the result of Strub {\rm(\cite{St})}, each isolated singular point $x$ 
 determines the equivalence class of a representation space $V_x$ 
 of a finite group $G_x$.
 \end{defi}
%
%%%%%%%%%%%%%%%%%%%%%%%
%
\bigskip
Let $x_1,\cdots ,x_k$ be 
the isolated singular points of $M$. 
 If $f \in {\cal D}(M)$, then $f(x_i) = x_i$ for $i=1,...,k$. 
Let $(U_{x_i}, \phi_{x_i})$ be a local chart of $M$ around $x_i$. 
There exists a smooth $\Gamma_{x_i}$-action on an open set 
$\tilde{U}_{x_i}$ in $V_{x_i}$ such that $\phi_{x_i}:\ U_{x_i} \to 
\tilde{U}_{x_i}/\Gamma_{x_i}$ is diffeomorphic. 
Let $C_i$ be the center of the group $\Gamma_{x_i}$. 
Then $C_i$ is naturally regarded as a subgroup of 
$Aut_{\Gamma_{x_i}}({V}_{x_i})$. 
Let $\overline{Aut}_{\Gamma_{x_i}}({V}_{x_i}) 
= Aut_{\Gamma_{x_i}}({V}_{x_i})/C_i.$
Take $\tilde{f}_i \in {\cal D}_{\Gamma_i}(\tilde{U}_i)$ such that 
$(\phi_{x_i}^{-1} \circ \pi_{x_i})_*(\tilde{f}_i) = f$  
on a neighborhood of $x_i$. 
Let 
$$\Psi : {\cal D}(M) \to \overline{Aut}_{\Gamma_{x_1}}({V}_{x_1})_0 
\times \cdots 
\times \overline{Aut}_{\Gamma_{x_k}}({V}_{x_k})_0$$
 be a homomorphism defined by  
$$\Psi(f)=(d\tilde{f}_1(\tilde{x}_1) C_1,\cdots ,
d\tilde{f}_k(\tilde{x}_k)C_k),$$
where $\tilde{x}_i$ is the isolated fixed point of 
$\tilde{U}_{x_i}$ with $\pi_i(\tilde{x}_i) = \phi_{x_i}(x_i)$. 
If $\dim M > 1$, using Satake \cite{Sa}, Lemma 1, we see that 
the homomorphism $\Psi$ is well-defined. If $\dim M =1$, 
it is clear since $\Gamma_{x_i} = \Z_2$.  
It is easy to see that $\Psi$ is epimorphic. 
%
%%%%%%%%%%% Lem 4.3 %%%%%%%%%%%%%%
%
\begin{lem}\label{lem4.3} 
Let $\varphi \in {\cal D}(U_{x_i})$ such that $\Psi(\varphi) =0$. 
Then $\varphi \in [\ker \Psi , {\cal D}(M)]$.
\end{lem}
%
%%%%%%%%%%%%%%%%
%
{\bf Proof}. \quad  
Let $\tilde{\varphi} \in {\cal D}_{\Gamma_{x_i}}(\tilde{U}_{x_i})$ 
such that $(\phi_{x_i}^{-1} \circ \pi_{x_i})_*(\tilde{\varphi}) = \varphi$.  
Since $\Psi(\varphi) =0$, $g=d\tilde{\varphi}(\tilde{x}_i) \in C_i$. 
There exists $\psi \in {\cal D}_{\Gamma_{x_i}}(\tilde{U}_{x_i})$ 
such that $\psi = g$ on a $\Gamma_{x_i}$ invariant 
neighborhood $\tilde{U}_{x_i}^0$ of $\tilde{x}_i$ in $\tilde{U}_{x_i}$. 
By Proposition \ref{prop3.4},  $\tilde{\varphi} \circ \psi^{-1} \in 
[\ker \Phi, 
{\cal D}_{\Gamma_{x_i}}(\tilde{U}_{x_i})]$. 
Since $(\phi_{x_i}^{-1} \circ \pi_{x_i})_*(\tilde{\varphi}\circ \psi^{-1}) 
= \varphi$ 
on a neighborhood $U_{x_i}^0 $ of $x_i$, $\varphi$ coincides with 
an element $\phi \in [\ker \Psi, {\cal D}(M)]$. 
Since $supp(\varphi \circ \phi^{-1}) \subset U_{x_i} 
\setminus U_{x_i}^0$, using Theorem \ref{th2.1}, (1), we can prove 
$\varphi \circ \phi^{-1} \in [\ker \Psi, {\cal D}(M)]$. 
This completes the proof of Lemma \ref{lem4.3}. 

\bigskip

 %
%%%%%%%% Th4.4 %%%%%%%%
%
\noindent
\begin{theo}\label{th4.4} \quad 
{\it Let $M$ be a smooth orbifold and let $x_1,\cdots ,x_k$ be the set of 
isolated singular points of $M$. Then the induced map  
$$\Psi_*: \ H_1({\cal D}(M)) \to 
\H_1(\overline{Aut}_{\Gamma_{x_1}}({V}_{x_1})_0)
\oplus\cdots 
\oplus H_1(\overline{Aut}_{\Gamma_{x_k}}({V}_{x_k})_0)$$}
is an isomorphism.
\end{theo}
%
%%%%%%%%%%%%%%%
%

\noindent
{\bf Proof}. \quad 
 Let  $f \in \ker(\Psi)$.  Then $f$ can be written as 
 $f = h_n \circ \cdots \circ h_1$ 
 with $h_j \in {\cal D}(U_{i_j}) \ (j=1,...,n)$ 
 such that $(U_{i_j},\ \phi_{i_j})$ is a local chart of $M$, 
 and $h_j \in ker \Psi$. 
 Using Theorem \ref{th2.1} (1) and Lemma \ref{lem4.3}, 
 $h_j \in [{\cal D}(M), {\cal D}(M)]$.
 Thus $\ker \Psi \subset [{\cal D}(M), {\cal D}(M)]$. 
 
 Since there is a short exact sequence
$$\begin{CD} 
\hspace{-1cm} 1 \to \ker \Psi \stackrel{\iota}{\to} {\cal D}(M)  
\stackrel{\Psi}{\to}  \overline{Aut}_{\Gamma_{x_1}}({V}_{x_1})_0 
\times \cdots \times 
\overline{Aut}_{\Gamma_{x_k}}({V}_{x_k})_0 \to 1 \end{CD},$$ 
we have the following exact sequence of homology groups: 

%\pagebreak

$$\begin{CD}
\hspace{-3cm} \ker \Psi/[\ker \Psi, {\cal D}(M)] \stackrel{\iota_*}{\to} 
H_1({\cal D}(M)) 
\end{CD}$$
$$\begin{CD}
\to \H_1(\overline{Aut}_{\Gamma_{x_1}}({V}_{x_1})_0)
\oplus\cdots 
\oplus H_1(\overline{Aut}_{\Gamma_{x_k}}({V}_{x_k})_0) \to 1.
\end{CD}$$
From the above argument we have $\iota_* = 0$, which completes the proof of 
Theorem \ref{th4.4}. 

\vspace{2mm}

%
%%%%%%%%%%% Remark 4.5 %%%%%%%%%%
%
\begin{rem} \quad 
Note that 
$$C_i/[C_i,\, Aut_{\Gamma_{x_i}}({V}_{x_i})_0] \to 
H_1(Aut_{\Gamma_{x_k}}({V}_{x_k})_0) \to 
H_1(\overline{Aut}_{\Gamma_{x_i}}({V}_{x_i})_0) \to 1
$$
is exact. Since $C_i$ is a finite group, by Corollary {\rm \ref{cor2.2}},
$H_1(\overline{Aut}_{\Gamma_{x_i}}({V}_{x_i})_0)$ is isomorphic to 
$H_1(Aut_{\Gamma_{x_k}}({V}_{x_k})_0)$. 
\end{rem}
%
%%%%%%%%%%%%%%%%
%

  Combining the fragmentation lemma (see  \cite {A-F1}, Lemma 1) 
  and Theorem \ref{th2.1} we have the following.  

 %
 %%%%%%%% Th 4.6 %%%%%%%%
 %
 \begin{theo} \label{th4.6}  
  Let $G$ be a finite group and 
 $M$ a smooth $G$-manifold. 
 If the orbit space $M/G$ has $\{G \cdot p_1,...,G \cdot p_k\}$ 
 as the isolated  singular points, then  
 $$H_1({\cal D}_G(M)) \cong   
  H_1(Aut_{G_{{p}_1}}(T_{{p}_1}M)_0) \times \cdots 
 \times  H_1(Aut_{G_{{p}_k}}(T_{{p}_k}M)_0).$$ 
 \end{theo}
  %
 %%%%%%%%%%%%%
 %

 Here we consider a smooth orbifold $M$ of dimension 2. In this case, 
 for each isolated singular point $x \in M$, the associate group $\Gamma_x$ 
 is isomorphic to a finite subgroup of $O(2)$.  
Combining Example \ref{ex2.4} with Theorem \ref{th4.4}, 
we have the following.

%
%%%%%%%%%%%%% Cor 4.7 %%%%%%%%%%%
%
\begin{cor} \quad  Suppose that $M$ has $n_1$ isolated singular points 
with associate group ${\bf Z}_2$, $n_2$ isolated singular points 
with associate group 
${\bf Z}_p (p \ge 3)$,  $n_3$ isolated singular points 
associated with ${\bf D}_2$ and 
$n_4$ isolated singular points with associate goroup ${\bf D}_r(r \ge 3)$. Then we have 
$$H_1({\cal D} (M)) \cong \overbrace{{\bf R} \times \cdots \times 
{\bf R}}^{n_1+n_2+2 n_3+n_4}
\times \overbrace{S^1 \times \cdots \times S^1}^{n_2}.$$ 
\end{cor}
%
%%%%%%%%%%%%%%%%
%

\vspace{5mm}

%\newpage
\noindent
{\bf 5. Application to compact foliations.}

\bigskip 
Let $M$ be an $n$-dimensional compact connected smooth manifold without boundary and ${\cal F}$ a codimension $q$ smooth foliation of $M$. A smooth diffeomorphism $f:M \to M$ is called {\it a foliation preserving diffeomorphism} (resp. {\it a leaf preserving diffeomorphism}) if for each point $x$ of $M$, the leaf 
through $x$ is mapped into the leaf through $f(x)$ (resp. $x$), that is, $f(L_x)=L_{f(x)}$ (resp. $f(L_x)=L_x$), where $L_x$ is the leaf of ${\cal F}$ which contains $x$. By ${\cal D} (M,{\cal F}) (resp. \,{\cal D} _L(M,{\cal F}))$ we denote the group of all foliation preserving (resp. leaf preserving) smooth diffeomorphisms of $(M,{\cal F})$ which are isotopic to the identity through foliation 
preserving (resp. leaf preserving) diffeomorphisms.

Then Tsuboi proved the following by looking at the proofs of Hermann{\cite H} and Thurston \cite {Th1}.

\bigskip
\noindent
{\bf Theorem 5.1}(Tsuboi \cite{Ts}, Rybicki \cite{R})\,. \quad {\it ${\cal D} _L(M,{\cal F})$ is perfect.}

\bigskip
A foliation with all leaves compact is said to be a compact foliation. A compact foliation is said to be a compact Hausdorff foliation if the leaf space $M/{\cal F}$ is Hausdorff. 
In this section we study the first homology of ${\cal D} (M,{\cal F})$ for compact Hausdorff foliations. 
Then we have the following nice local picture for a compact Hausdorff smooth foliation of codimension $q$. 

%%%%%%%%%%%%%%
\bigskip
\noindent
{\bf Proposition 5.2}(Epstein \cite {E})\,. \quad {\it There is a generic leaf $L_0$ with property that there is an open dense subset of $M$, where the leaves have all trivial holonomy and are all diffeomorphic to $L_0$. Given a leaf $L$, we can describe a neighborhood $U(L)$ of $L$, together with the foliation on the neighborhood as follows. There is a finite subgroup $G(L)$ of $O(q)$ such that $G(L)$ acts freely on $L_0$ on the right and $L_0/G(L) \cong L$. Let $D^q$ be the unit disk. We foliate $L_0 \times D^q$ with leaves of the form $L_0 \times \{ pt \}$. This foliation is preserved by the diagonal action of $G(L)$, defined by $g(x,y)=(x \cdot g^{-1},g \cdot y)$ for $g \in G(L)$, $x \in L_0$ and $y \in D^q$. So we have a foliation induced on $U=L_0 \times_{G(L)} D^q$. The leaf corresponding to $y=0$ is $L_0/G(L)$. Then there is an embedding $\varphi : U \rightarrow M$ with $\varphi (U)=U(L)$, which preserves leaves and $\varphi (L_0/G(L))=L$.}

\bigskip
\noindent
{\bf Definition 5.3}\,. \quad A leaf $L$ in ${\cal F}$ is called a {\it singular leaf} if $G(L)$ is not trivial. A singular leaf  $L$ is called an {\it isolated singular leaf} if the origin of $D^q$ is an isolated fixed point of $G(L)$.

\bigskip
%%%%%%%%%%%%%%%%%%%%%%%
Note that there are finitely many isolated compact leaves in ${\cal F}$ since $M$ is compact.
Then we have the following. 

\bigskip
\noindent
{\bf Theorem 5.4}\,. \quad {\it Let ${\cal F}$ be a compact Hausdorff foliation of $M$ and $L_1, \cdots ,L_k$ be all isolated singular leaves of ${\cal F}$. Then we have 
$$H_1({\cal D} (M,{\cal F}))\cong H_1(Aut_{G(L_1)}(D^q)_0)\oplus \cdots \oplus H_1(Aut_{G(L_k)}(D^q)_0).$$}
%where $Z(G(L_i))=\{ A\in GL_+(q,{\bf R}) \mid g\cdot A=A\cdot g \, (\forall g \in G(L_i)\}$.}

\bigskip
We consider the case $q=1$. Then we have the following. 
  
\bigskip
\noindent
{\bf Corollary 5.5}\,.\quad {\it Let ${\cal F}$ be a codimension one compact foliation of $M$. 

\noindent
$(1)$  If ${\cal F}$ is transversely orientable, then ${\cal D} (M,{\cal F})$ is perfect.

\noindent
$(2)$  If ${\cal F}$ is not transversely orientable, then $H_1({\cal D} (M,{\cal F})) \cong {\bf R}\times {\bf R}$.}

\bigskip
\noindent
{\bf Proof}.\, \quad For the case ${\cal F}$ is transversely orientable, ${\cal F}$ is a bundle foliation. In this case the proof follows from the first part of the proof of Theorem 5.4 below. Next we consider the case ${\cal F}$ is not transversely orientable. 
A non-trivial finite subgroup of $O(1)$ is $O(1)$ itself, which is generated by the reflection. Then we have $Aut_{O(1)}(D^1)_0=GL_+(1,{\bf R})$. Since ${\cal F}$ has only two singular leaves in this case, the proof follows from Theorem 5.4.

\bigskip
We consider the case $q=2$. Then finite subgroups $G$ of $O(2)$ are as follows: 

$G$ is either a group of $k$ rotations which is isomorphic to ${\bf Z}_k$ $(k\ge 2)$ or a group of $l$ rotations and $l$ reflections which is isomorphic to ${\bf D}_l=\{ u,v \,; u^l=v^2=(uv)^2=1 \}$ $(\ell \ge 1)$. 

\bigskip
\noindent
{\bf Definition 5.6}\,. \quad A singular leaf $L$ in ${\cal F}$ is called a rotation leaf with holonomy ${\bf Z}_p$ or a dihedral leaf with holonomy  ${\bf D}_l$ according that $G(L)$ is isomorphic to ${\bf Z}_k$ or ${\bf D}_l (\ell \ge 2)$. 

\bigskip
Then we have the following. 
  
\bigskip
\noindent
{\bf Corollary 5.7}\,.\quad {\it Let ${\cal F}$ be a codimension two compact foliation of $M$. Suppose that ${\cal F}$ has $n_1$ rotation leaves with holonomy ${\bf Z}_2$, $n_2$ rotation leaves with holonomy ${\bf Z}_p (p \ge 3)$, $n_3$ dihedral leaves with holonomy ${\bf D}_2$ and $n_4$ dihedral leaves with holonomy ${\bf D}_{\ell} (\ell \ge 3)$. Then we have 

$H_1({\cal D} (M,{\cal F})) \cong \overbrace{{\bf R} \times \cdots \times {\bf R}}^{n_1+n_2+2n_3+n_4}\times \overbrace{ S^1 \times \cdots \times S^1}^{n_2}.$
}

\bigskip
\noindent
{\bf Proof}.\, \quad This follows from Theorem 5.4 and Example 2.4.

\vspace{5mm}

%\newpage
\noindent
{\bf 6. Proof of Theorem 5.4}

\bigskip

 First we consider the case where ${\cal F}$ has no singular leaves. Then we shall prove that $H_1({\cal D} (M,{\cal F}))$ vanishes, that is, ${\cal D} (M,{\cal F})$ is perfect. In this case, we have a smooth fiber bundle $p:M \to M/{\cal F} \cong B$, where $B$ ia a compact smooth manifold of dimension $n-q$. We denote by ${\cal D} (B)$ the group of all smooth diffeomorphisms of $B$ which are isotopic to the identity. Then we define a map $p_*:{\cal D} (M,{\cal F}) \to {\cal D} (B)$ by $p_*(f)(\bar{x})=p\circ f(x)$ for 
any $f\in {\cal D} (M,{\cal F})$ and $\bar{x} \in B$, $x\in M$ with $p(x)=\bar{x}$.

Then we have the following Lemma, which is easily proved and is left to the reader.

\bigskip
\noindent
{\bf Lemma 6.1}\,. \quad {\it $p_*$ is an epimorphism.}

\bigskip
Take any $f\in {\cal D} (M,{\cal F})$. Then we may assume $f$ is close to the identity $1_M$. Thus $\bar{f}=p_*(f)$ is close to the identity $1_B$. From the perfectness of ${\cal D} (B)$ (Thurston \cite {Th1}), there exist 
$\bar{g}_1,\bar{g}_2,\cdots ,\bar{g}_{2k} \in {\cal D} (B)$ satisfying $\bar{f}=\displaystyle\prod_{i=1}^s [\bar{g}_{2i-1}, \bar{g}_{2i}]$. 
 From Lemma 6.1, there are $g_1,g_2, \cdots ,g_{2k} \in {\cal D} (M,{\cal F})$ satisfying that $p_*(g_i)=\bar{g}_i \, (i=1,2,\cdots ,2k)$. Then since $p_*(f\circ (\displaystyle\prod_{i=1}^s [g_{2i-1}, g_{2i}])^{-1})=1_B$ and $f\circ (\displaystyle\prod_{i=1}^s [g_{2i-1}, g_{2i}])^{-1}$ is close to the identity $1_M$,  $f\circ (\displaystyle\prod_{i=1}^s [g_{2i-1}, g_{2i}])^{-1}$ is contained in the connected component of $\ker p_*$, which is ${\cal D} _L(M,{\cal F})$. From Theorem 5.1, $f\circ (\displaystyle\prod_{i=1}^s [g_{2i-1}, g_{2i}])^{-1}$ is expressed as a composition of commutators of elements in ${\cal D} _L(M,{\cal F})$. This proves that ${\cal D} (M,{\cal F})$ is perfect.

%%%%%%%%%%%%%%%%%%%%%
\bigskip
Next we shall prove the case where ${\cal F}$ has isolated singular leaves. Let $L$ be an isolated singular leaf and $U(L)$ a saturated neighborhood of $L$, which is foliation preserving diffeomorphic to $U=L_0 \times_{G(L)} D^q$ as in Proposition 5.2. We identify $U(L)$ with $U$. 

Take $f \in {\cal D} (M,{\cal F})$. Note that $f$ preserves $L$. Then $f\mid_U$ can be lifted to a foliation preserving and $G(L)$-equivariant embedding $\tilde{f} : L_0 \times  D^q \to L_0 \times  D^q(r)$ up to $G(L)$-action for some $r>0$, where $D^q(r)$ is the disk of radius $r$, centered at the origin. Thus we have a $G(L)$-equivariant embedding $\bar{f} : D^q \to D^q(2)$. We assign $f$ to the differential of $\bar{f}$ at the origin, $d\bar{f}(0)$. We denote this map by $\Phi_L :{\cal D} (M,{\cal F}) \to Aut_{G(L)}(D^q)_0 $ which is a homomorphism. 

\bigskip
\noindent
{\bf Lemma 6.2}\,. \quad {\it The map $\Phi_L $ is surjective.}

\bigskip
\noindent
{\bf Proof}\,. \quad For any $A \in Aut_{G(L)}(D^q)_0$, we can get a $G(L)$-equivariant diffeomorphism $\bar{f} : V^q \to V^q$ satisfying that $d\bar{f}(0)=A$ and $\bar{f}$ is supported in $D^q$, where $V^q$ is the representation space of $G(L)$. Then we define a foliation preserving and $G(L)$-equivariant diffeomorphism $\tilde{f} : L_0 \times  D^q \to L_0 \times  D^q$ by $\tilde{f}(x,y)=(x,\bar{f}(y))$, where $(x,y)$ denotes the local coordinate of $L_0 \times  D^q$. $\tilde{f}$
induces the foliation preserving diffeomorphism $f\mid_U : L_0 \times_{G(L)}  D^q \to L_0 \times_{G(L)}  D^q$. Then we can extend $f\mid_U$ to a foliation preserving diffeomorphism of $(M,{\cal F})$ by putting the identity outside of $U$. This completes the proof. 

\bigskip

\noindent
{\bf Proposition 6.3}\,. \quad {\it Let $L$ be an isolated singular leaf and $U$ a saturated neighborhood of $L$. Suppose that ${f \in \ker \Phi_{L} }$ is supported in $U$. Then $f$ can be expressed as a composition of commutators of elements with support in $U$ in $\ker \Phi_{L}$ and ${\cal D} (M,{\cal F})$, that is, $f \in [\ker \Phi_{L}, {\cal D} (M,{\cal F})]$.}

\bigskip
\noindent
{\bf Proof}\,. \quad Let ${f \in \ker \Phi_L }$. Then from Proposition 3.4, the $G(L)$-equivariant diffeomorphism $\bar{f} : D^q \to D^q$ can be expressed as a composition of commutators of elements which are extendable to elements in $\ker \Phi$ and ${\cal D} _{G(L_i)}(V^q)$, say $\bar{f}=\displaystyle\prod_{j=1}^s [\bar{f}_{2j-1}, \bar{f}_{2j}]$, where $\bar{f}_{2j-1} \in \ker \Phi$ and $\bar{f}_{2j} \in {\cal D} _{G(L_i)}(V^q)$. From Lemma 6.2, we can lift $\bar{f}_j\,(j=1,2,\cdots ,2s)$ to foliation preserving diffeomorphisms $f_j\,(j=1,2,\cdots ,2s)$ satisfying that each $f_j$ is supported in $U$ and $f_{2\ell-1} \in \ker \Phi_{L} (1 \le \ell \le s)$. This completes the proof. 

\bigskip
\noindent
{\bf Proposition 6.4}\,. \quad {\it Let $L$ be a non-isolated singular leaf and $U$ a saturated neighborhood of $L$. Suppose that ${f \in \ker {\cal D} (M,{\cal F})}$ is supported in $U$. Then $f$ can be expressed as a composition of commutators of elements with support in $U$ in ${\cal D} (M,{\cal F})$.}

\bigskip
\noindent
{\bf Proof}\,. \quad Let $\bar{f} : D^q \to D^q$ be the $G(L)$-equivariant diffeomorphism induced from $f$ as above. Note that $\dim V^{G(L)} > 0$. From Theorem 2.1(2), we see that $\bar{f} : D^q \to D^q$ can be expressed as a composition of commutators of elements in ${\cal D} _{G(L)}(V^q)$ with support in $D^q$. As in the proof of Proposition 6.2, these elements can be lifted to foliation preserving diffeomorphisms of $(M,{\cal F})$ with support in $U$. This completes the proof. 

\bigskip
\noindent
{\bf Proof of Theorem 5.4}\,. \quad  Let $L_1,\cdots ,L_k$ be all isolated singular leaves of ${\cal F}$. Let $$\Phi : {\cal D} (M,{\cal F}) \to Aut_{G(L_1)}(D^q)_0\times \cdots \times  Aut_{G(L_k)}(D^q)_0$$ be the homomorphism defined by 
$\Phi(f)=(\Phi_{L_1}(f),\cdots ,\Phi_{L_k}(f))$ for any $f \in {\cal D} (M,{\cal F})$. 
Take $f \in \ker \Phi $. Then by using the fragmentation argument, there are $f_1, \cdots , f_n \,\in \ker \Phi \,(n > k)$ such that 

\noindent
$(1)$ $f_j$ is supported in a saturated neighborhood $U_j$ of the isolated singular leaf $L_j$ for $1 \le j \le k$ and is supported in a saturated neighborhood $U_j$ of a non-isolated singular leaf or non-singular leaf $L_j$ for $j > k$, and

\noindent
$(2)$ $f=f_n \circ \cdots \circ f_1$.

Since there is a short exact sequence
$$\begin{CD} 
\hspace{-1cm} 1 \to \ker \Phi \stackrel{\iota}{\to} {\cal D}(M,{\cal F})  
\stackrel{\Phi}{\to}   Aut_{G(L_1)}(D^q)_0 \times \cdots \times  Aut_{G(L_k)}(D^q)_0 \to 1 \end{CD},$$ 
we have the following exact sequence of homology groups: 
$$\begin{CD}
\hspace{-3cm} \ker \Phi/[\ker \Phi, {\cal D}(M,{\cal F})] \stackrel{\iota_*}{\to} 
H_1({\cal D}(M,{\cal F})) 
\end{CD}$$
$$\begin{CD}
\to H_1(Aut_{G(L_1)}(D^q)_0 \times \cdots \times 
Aut_{G(L_k)}(D^q)_0) \to 1.
\end{CD}$$
Then since the fact of the non-singular case, Proposition 6.3 and Proposition 6.4 imply the equality $\ker \Phi=[\ker \Phi, {\cal D}(M,{\cal F})]$, the proof is complete.

\bigskip

\noindent
{\bf 7. Application to an $S^1$-action on $S^3$}

\bigskip

In this section we shall apply Theorem \ref{th2.1} to the case of a locally free $S^1$-action on the 3-dimensional sphere.

We consider $S^1$ as the set of unit complex numbers $U(1)$ and put 
$S^3=\{(w_1,w_2) \in {\bf C} \mid |w_1|^2+|w_2|^2=1 \}$ with $U(1)$-action given by 
$$z\cdot (w_1,w_2)=(zw_1,z^2w_2), \quad z\in U(1).$$ 
Then it has two orbit types $\{(1),\, {\bf Z}_2 \}$ and the orbit space $S^3/U(1)$ is homeomorphic to the space known as the tear drop which is the two dimensional sphere with one isolated singular point. Let $V^1$ denote the non-trivial one dimensional representation space of ${\bf Z}_2$. Note that the slice representation at the point $e=(0,1) \in S^3$ coincides with $V \oplus V$. 
Then by Example \ref{ex2.3} we have

\bigskip
\noindent
{\bf Proposition 7.1}\,. \quad {\it $H_1({\cal D}(S^3/U(1))) \cong {\bf R}$.
}

\bigskip
In the rest of this section we shall prove the following.

\bigskip
\noindent
{\bf Theorem 7.2}\,. \quad {\it $H_1({\cal D}_{U(1)}(S^3)) \cong {\bf R}\times U(1)$.
}

\bigskip
  Let $W$ be an invariant tubular neighborhood of the exceptional orbit 
$U(1) \cdot e$.  Put $\tilde{D} = \{(x,y) \in {\bf R}^2 | \ x^2 + y^2 < 1\}$.Then $W$ can be identified with $\tilde{D} \times S^1$ 
on which the $U(1)$-action is given by 
$$
z \cdot (w_1, w_2) = (zw_1, z^2 w_2), \qquad  
z \in U(1),\ (w_1, w_2) \in \tilde{D} \times S^1.
$$ 
We can consider the group ${\cal D}_{U(1)}(W)$ and ${\cal D}(W/U(1))$ as the subgroups of ${\cal D}_{U(1)}(S^3)$ and ${\cal D}(S^3/U(1))$ respectively.
Let 
  $$P:\ {\cal D}_{U(1)}(S^3) \to {\cal D}(S^3/U(1))$$
be the canonical map induced from the natural projection from $S^3$ to 
$S^3/U(1)$.  \par 
   Note that  $\{x^2, y^2, xy\}$ is a Hilbert basis of 
${\bf R}[V^2]^{{\bf Z}_2}$.  Let $\pi:\ W \to {\bf R}^3$ be the map defined by 
$$
  \pi(x,y,w_2) = (x^2,y^2,xy), \qquad  \ (x,y,w_2) \in W.  
$$
Let $\iota: \ \tilde{D}/{\bf Z}_2 \to W/U(1)$ be the map given by 
$$
\iota({\bf Z}_2 \cdot w) = U(1) \cdot (w,1), \qquad  \ 
w \in \tilde{D}.
$$  

%%%%%%%%%%%%%
\bigskip
Here we refer the smooth structure of the orbit space of a 
representation space of a Lie group which was studied by 
Bierstone \cite{Bi1} and Schwarz \cite{Sc1}, \cite{Sc2}. 
Note that $\tilde{D}$ is a slice of the orbit $U(1) \cdot e$.  
Then the inclusion map $i:\ \tilde{D} \hookrightarrow W$ induces the 
isomorphism 
$i^* :\ C^{\infty}_{U(1)}(W) \to C^{\infty}_{{Z}_2}(\tilde{D})$ 
of the algebras of the invariant functions spaces.  
Then by the definition of the smooth structure of the orbit space, the map 
$\iota$ is diffeomorphic. 
Let $p: \ \tilde{D} \to {\bf R}^3$ denote the map given by 
$p(x,y) = (x^2,y^2,xy).$  By Bierstone \cite{Bi2}, $\S 2$, 
$p(\tilde{D})$ is 
diffeomorphic to the orbit space $\tilde{D}/{\bf Z}_2.$   Therefore we have 
the diffeomorphisms of the spaces 
$$ W/U(1) \cong \tilde{D}/{\bf Z}_2 \cong p(\tilde{D}). $$
 
Let $P_0 :\ {\cal D}_{{\bf Z}_2}(\tilde{D}) \to {\cal D}(p(\tilde{D}))
  \cong {\cal D}(\tilde{D}/{\bf Z}_2)$ be the canonical map induced from $p$. 
   Here we consider the manifold 
$(\tilde{D} \times U(1))/{\bf Z}_2$ with $U(1)$-action given by 
$$
z \cdot [w,u] = [w, zu], \qquad  \ z \in U(1),\ 
(w,u) \in \tilde{D} \times U(1)
$$ 
where $ [w,u] $ denotes the equivalence class ${\bf Z}_2 \cdot (w,u).$
There exists an equivariant diffeomorphism  
$\gamma : \ (\tilde{D} \times U(1))/{\bf Z}_2 \to W$ defined by 
$\gamma ([w,u]) = (uw, u^2)$. We identify $U(1)$ with $S^1$ by the 
canonical way. Then $\gamma$ induces the isomorphism 
$$\gamma_*:\ {\cal D}_{U(1)}((\tilde{D} \times U(1))/{\bf Z}_2) \to {\cal D}_{U(1)}(W)$$
given by $\gamma_*(h) = \gamma \circ h \circ \gamma^{-1}$.  

\bigskip
  Let $q:\ (\tilde{D} \times U(1))/{\bf Z}_2 \to \tilde{D}/{\bf Z}_2$ denote 
the natural projection. Let $\{w\}$ denote the 
orbit ${\bf Z}_2 \cdot w$ of $w \in \tilde{D}$.
Set $\tilde{D}_0 = \tilde{D} \setminus \{0\}.$  
Then the restriction $q:\ (\tilde{D}_0 \times U(1))/{\bf Z}_2 \to 
\tilde{D}_0/{\bf Z}_2$ is a smooth principal $U(1)$-bundle. Let $s$ be a 
smooth section of the bundle $q$ defined by 
$$s(\{w\}) = [w,\ \frac{\bar{w}}{|w|}], \qquad  \ 
w \in \tilde{D}_0.$$ 
  Choose $h \in Ker P$. Then $h$ (resp. $\gamma_*^{-1}(h)$) preserves the orbits 
of $W$ (resp. $\tilde{D}_0 \times U(1))/{\bf Z}_2$). Thus there exists a 
smooth map $a_h:\ \tilde{D}_0/{\bf Z}_2 \to U(1)$ such that 
$$\gamma_*^{-1}(h)(s(\{w\})) =  a_h(\{w\})\cdot s(\{w\}), \qquad  
  w \in \tilde{D}_0.$$     
It follows that 
$$
\gamma_*^{-1}(h)([w,u]) = [w, a_h(\{w\})u], \qquad 
[w,u] \in \tilde{D}_0 \times U(1)/{\bf Z}_2.
$$  

\bigskip
\noindent
{\bf Lemma 7.3}\,. \quad {\it The map $a_h$ can be extended to a smooth map $\bar{a}_h$ over $\tilde{D}/{\bf Z}_2$.}

\bigskip
\noindent
{\bf Proof}\,. \quad If $w \in \tilde{D}_0,\ u \in S^1$, then 
\begin{eqnarray*}
 h(w,u^2) &=& (\gamma \circ \gamma_*^{-1}(h) \circ \gamma^{-1})(w, u^2) \\
           &=& (\gamma \circ \gamma_*^{-1}(h))([\bar{u} w, u]) \\
           &=& \gamma([\bar{u}w, \, a_h(\{\bar{u} w\})u]) \\
           &=& (a_h(\{\bar{u} w\})w, \ a_h(\{\bar{u} w\})^2 u^2).
\end{eqnarray*}
Let $\rho_1:\ W \to \tilde{D}$ be the projection on the first factor and 
let $h_1:\ W \to \tilde{D}$ be the composition $\rho_1 \circ h$.  Since 
$h$ preserves the orbit $U(1) \cdot e$, 
$h_1(0,u) = 0$ for any $u \in S^1$.  
Let $\hat{a}_h:\ \tilde{D}_0 \to U(1)$ be a map defined by 
$\hat{a}_h(w) = a_h(\{w\})$.  Let 
$\phi:\ \tilde{D} \to \tilde{D}$ be a map defined by $\phi(w) 
= h_1(w,1)$.  Then $\phi$ is a smooth map satisfying  
  $$\phi(w) = \left\{ \begin{array}{lll}
                  \hat{a}_h(w)w & & {\rm if} \ w \neq 0,\\
                  0   & & {\rm if} \ w = 0. 
                  \end{array}  \right.$$  
By the Taylor's formula we have 
$$ \phi(w) = \left( \int_0^1 \phi'(tw) \, dt \right) w, \qquad  \ 
w \in \tilde{D}. $$ 
It follows that 
    $$(7.1) \qquad \qquad \hat{a}_h(w) = \int_0^1 \phi'(tw) \, dt,   
 \qquad  \ w \in \tilde{D}_0.  $$ 
Thus the map $\hat{a}_h$ is extended to a smooth map $\tilde{a}_h$ over 
$\tilde{D}$  defined by the right hand side of the equation (7.1).
 Let 
$$ 
\bar{a}_h(\{w\}) = \tilde{a}_h(w), \qquad  \ w \in \tilde{D}. 
$$ 
By definition $\bar{a}_h$ is a smooth extension of $a_h$ over 
$\tilde{D}/{\bf Z}_2$.  This completes the proof of Lemma 7.3.

\bigskip
 Let $h \in {\cal D}_{U(1)}(W)$.  Since $h$ preserves the 
exceptional orbit $U(1) e$, there exists $z_h \in U(1)$ satisfying 
$h(e) = z_h \cdot e$. Since the isotropy subgroup of $U(1)$ at the point $e$ 
is ${\bf Z}_2= \{\pm 1\}$, we have a map 
$$ T:\ {\cal D}_{U(1)}(W) \to U(1)/{\bf Z}_2 $$ 
defined by $T(h) = z_h {\bf Z}_2 $.  \par 

%%%%%%%delete
%  For $h \in D_{U(1)}(W)$ let $\tilde{h} = z_h^{-1} \cdot h$. 
%Then $\tilde{h}(e) = e$ and the differential $(d\tilde{h})_e$ of $\tilde{h}$ at the %point $e$ induces the ${\Z}_2$-invariant linear automorphism 
%$\varphi_h$ of the normal space of the orbit $U(1)\cdot e$ at $e$.  
%Since the normal space is isomorphic to ${\R}^2$ as a representation 
%space ${\Z}_2$, we can define the homomorphism 
%  $$ \Phi: \  D_{U(1)}(S^3) \to GL({\tR}^2), \quad 
%  \Phi(h) = \varphi_h {\Z}_2.$$ 
  
\bigskip
\noindent
{\bf Lemma 7.4}\,. \quad {\it 
$(1)$ \ The map $T$ is a group homomorphism. \par 
$(2)$ \ If $h \in Ker P$, then $T(h) = \bar{a}_h(0) \ {\bf Z}_2$.}

\bigskip
\noindent
{\bf Proof}\,. \quad  
(1) \ For $h_1, h_2 \in {\cal D}_{U(1)}(W)$, we have 
$$
 z_{h_1 \circ h_2} \cdot e = (h_1 \circ h_2)(e) =  z_{h_1} \cdot (z_{h_2} \cdot e) 
 = (z_{h_1} z_{h_2}) \cdot e. 
$$
Thus $T$ is a group homomorphism. \par 
  (2) \ Let $h \in Ker P$.  As in the proof of Lemma 7.3, if 
 $w \in \tilde{D}_0,\ u \in S^1$, then  
 $$ h(w,u^2) = (\bar{a}_h(\{\bar{u} w\})w, \ \bar{a}_h(\{\bar{u} w\})^2 u^2). $$ 
 Thus we have 
 $$ T(h) \cdot e = h(0, 1) = (0,\ \bar{a}(0)^2) = \bar{a}(0) \cdot e. $$ 
 This completes the proof of Lemma 7.4.

\bigskip
 
%%%%%%%%%%%%%%%%%%%% \eta %%%%%%%%%%%%%%%%%%%%%%%%%
%   Note that $-1 \in U(1)$ induces an element 
% $\eta \in D_{U(1)}(W)$ given by 
%$$\eta((w,u)) =(-1) \cdot (w,u) = (-w,u) \quad {\rm for}\ (w,u) \in W.$$ 
%Then $\{\pm \eta\}$ is a central subgroup of $ D_{U(1)}(W)$. 
%%%%%%%%%%%%%%%%%%%%%%%%%%%%%%%%%%%%%%%%%%%%%%%%% 
 Let 
 $$\Theta :\ {\cal D}_{U(1)}(W) \to {\cal D}(W/U(1)) \times 
 U(1)/{\bf Z}_2$$ 
 be a map defined by $\Theta(h) = (P(h), T(h)).$  
 
\bigskip
\noindent
{\bf Lemma 7.5}\,. \quad {\it 
The map $\Theta$ is a surjective group homomorphism.}

\bigskip
\noindent
{\bf Proof}\,. \quad 
 Let $f \in  {\cal D}(W/U(1)) \cong {\cal D}(\tilde{D}/{\bf Z}_2), 
  \ z \in U(1)$. There exists $\hat{f} \in {\cal D}_{{\bf Z}_2}(\tilde{D})$ 
  such that 
  $P_0(\hat{f}) =f$.  Define 
  $$ h([w,u]) = [\hat{f}(w),\, zu] \quad {\rm for} \ [w,u] \in 
  (\tilde{D} \times U(1))/{\bf Z}_2.  $$ 
  Then it is easy to see that 
  $\Theta(\gamma_*^{-1}(h) ) = (f,\,z \, {\bf Z}_2)$ 
  and Lemma 7.5 follows.   

\bigskip
\noindent
{\bf Lemma 7.6}\,. \quad {\it 
Let $a \in C^{\infty}(p(\tilde{D}))$ with $a(0)=0$. For 
 $\delta >0$, 
there exists $v \in C^{\infty}(p(\tilde{D}))$ and 
$f \in {\cal D}(p(\tilde{D}))$ 
such that \par 
$(1)$ \ $v \circ f - v = a$ on some neighborhood of $0$, \par 
$(2)$ \ $|v(\{w\})| \le \delta \quad (w \in \tilde{D}),$ \par 
$(3)$ \ $v$ has a compact support.}

\bigskip
\noindent
{\bf Proof}\,. \quad 
 Since $a(0)=0$, the function $a$ has of the form 
 $$a(X,Y,Z) = c_1 X + c_2 Y + c_3 Z + \xi (X,Y,Z),$$ 
 where $\xi$ is the higher order term. Set $\hat{a}=a \circ p$. Then 
 $\hat{a}$ is a smooth map and 
 $$\hat{a}(x,y) =  c_1 x^2 + c_2 y^2 + c_3 xy 
 + \xi (x^2,y^2,xy)  \quad {\rm for} \ x,y \in \tilde{D}.$$ 
 Suppose that the homogeneous terms of degree 2 of $\hat{a}$ is not 
 identically zero. Then there exists a linear transformation  
 $$ (7.2) \qquad 
 \left\{
 \begin{array}{lll}
     \tilde{x} & = & d_1 x + d_2 y \\
     \tilde{y} & = & d_3 x + d_4 y,
 \end{array}
 \right.
 $$
 such that  $\hat{a}(x,y)$ is written as 
 $$ \hat{a}(x,y) = \check{a}(\tilde{x},\tilde{y}) 
 = \alpha \tilde{x}^2 + \beta \tilde{y}^2 
 + \eta(\tilde{x}^2, \tilde{y}^2, \tilde{x} \tilde{y}),  $$
 where $\eta$ is the higher order term.  Then $\check{a}$ is written as 
 $$ \check{a}(\tilde{x},\tilde{y}) 
 = \tilde{x}^2 \mu(\tilde{x}^2, \tilde{y}^2, \tilde{x} \tilde{y}) 
 + \tilde{y}^2 \nu(\tilde{x}^2, \tilde{y}^2, \tilde{x} \tilde{y}), $$
 where $\mu$ and $\nu$ are the smooth functions.  Note that the 
 argument is holds if the homogeneous terms of degree 2 of $\hat{a}$ is 
 identically zero. \par 
   The linear transformation $(7.2)$ induces the following linear 
transformation 
  $$   (7.3) \qquad 
 \left\{
 \begin{array}{lll}
     \tilde{X} & = & d_1 X + d_2 Y + 2 d_1 d_2 Z,\\
     \tilde{Y} & = & d_3 X + d_4 Y + 2 d_3 d_4 Z, \\
     \tilde{Z} & = & d_1 d_3 X + d_2 d_4 Y +  (d_1 d_4+ d_2 d_3) Z,
 \end{array}
 \right.
 $$
 Then $a(X,Y,Z)$ is written as
 $$a(X,Y,Z) = \tilde{a}(\tilde{X}, \tilde{Y}, \tilde{Z}) 
 = \tilde{X} \mu(\tilde{X}, \tilde{Y}, \tilde{Z}) 
 + \tilde{Y} \nu(\tilde{X}, \tilde{Y}, \tilde{Z}).  $$
  Let $$h_0(\tilde{x}, \tilde{y}) = 
  (\tilde{x}\sqrt{1 + \mu(p(\tilde{x}, \tilde{y}))},
  \tilde{y}\sqrt{1 + \nu(p(\tilde{x}, \tilde{y}))}).$$ 
  Then $h_0$ is a local ${\bf Z}_2$-equivariant diffeomorphism on a 
  sufficiently small neighborhood of $0$ in $\tilde{D}$ which is 
  equivariantly isotopic to the identity.  Then there exists 
  $h \in {\cal D}_{{\bf Z}_2}(\tilde{D})$ which  coincides with $h_0$ on a small 
  neighborhood of $0$.  Set $\tilde{f} = P_0(h)  
   \in {\cal D}(p(\tilde{D})).$  Then  
 \begin{eqnarray*}
 \lefteqn{ \tilde{f}(\tilde{X}, \tilde{Y}, \tilde{Z})  =  
   \left( \tilde{X}(1+ \mu(\tilde{X}, \tilde{Y}, \tilde{Z})), \       
  \tilde{Y}(1+ \nu(\tilde{X}, \tilde{Y}, \tilde{Z})),\right.} \\  
 &\hspace{50mm} & \left. \tilde{Z} \sqrt{(1+ \mu(\tilde{X}, \tilde{Y}, \tilde{Z}))
  (1+ \nu(\tilde{X}, \tilde{Y}, \tilde{Z})}) \right) 
 \end{eqnarray*}
 on some small neighborhood of $0$.  
 Let $\tilde{v}:\ p(\tilde{D}) \to {\bf R}$ be a smooth function such that 
 \begin{description} \par 
 \item \quad (1) \   $\tilde{v}(\tilde{X}, \tilde{Y}, \tilde{Z}) 
 = \tilde{X} + \tilde{Y}$ on a neighborhood $U$ of $0$, \par 
 \item \quad (2) \  $|\tilde{v}(\{w\})| \le \delta \quad  
 (w \in \tilde{D})$, \par 
 \item \quad (3) \ $\tilde{v}$ has a compact support.
 \end{description}
 Then 
 $\tilde{v} \circ \tilde{f} - \tilde{v} = \tilde{a}$ on some neighborhood of 
 $0$ in $p(\tilde{D})$.  From the linear transformation (7.3), 
 we have the maps $v$ and $f$ such that 
 \begin{eqnarray*} 
 v(X,Y,Z) &=& \tilde{v}(\tilde{X}, \tilde{Y}, \tilde{Z}), \\ 
  f(X,Y,Z) &=& \tilde{f}(\tilde{X}, \tilde{Y}, \tilde{Z}).  
 \end{eqnarray*} 
 Then $v \circ f - v = a$ on some neighborhood of $0$ and Lemma 7.6 
 follows. 
 
\bigskip
\noindent
{\bf Proposition 7.7}\,. \quad {\it 
$ Ker \, \Theta$ is contained in $[{\cal D}_{U(1)}(W),\ {\cal D}_{U(1)}(W)]$.}

\bigskip
\noindent
{\bf Proof}\,. \quad 
  Let $h \in {\cal D}_{U(1)}(W)$ with $\Theta(h)= 0$. By the fragmentation lemma (\cite{A-F1}, Lemma 1), we can assume that $h=h_1\circ h_2$ such that $h_1\in {\cal D}_{U(1)}(W)$ and $h_2\in {\cal D}_{U(1)}(S^3-U(1)\cdot e)$. Since $U(1)$ acts freely on $S^3-U(1)\cdot e$, by the main theorem in \cite{A-F1}, $h_2\in [{\cal D}_{U(1)}(S^3-U(1)\cdot e), {\cal D}_{U(1)}(S^3-U(1)\cdot e)]$. Thus we can assume $h\in {\cal D}_{U(1)}(W)$ with $\Theta(h)=0$. Note that $ h\in \, Ker \, P$ and 
$T(h) = {\bf Z}_2$.  By Lemma 7.4 (2)
$\bar{a}_h(0) = 1$ or $-1$.

First we assume that $\bar{a}_h(0) = 1.$  
We shall identify $\tilde{D}/{\bf Z}_2$ with $p(\tilde{D})$.  
  Let $\hat{a}_h:\ p(\tilde{D}) \to {\bf R}$ be 
a smooth function such 
that $\exp (2 \pi \hat{a}_h(x)\sqrt{-1}) 
= \bar{a}_h(x) \ (x \in p(\tilde{D}))$ and 
$\hat{a}_h(0)=0.$    
By Lemma 7.6 there exist $f \in {\cal D}(p(\tilde{D}))$ and 
$v \in C^{\infty}(p(\tilde{D}))$ such that \par 
\begin{description}
\item $(1)$ \ $v \circ f - v = \hat{a}$ on some neighborhood $U$ of $0$, 
\par 
\item $(2)$ \ $|v(p(w))| \le \frac{1}{2} \quad (w \in \tilde{D}),$ \par 
\item $(3)$ \ $v$ has a compact support. 
\end{description} \par 
  Let $\hat{f} \in  {\cal D}_{{\bf Z}_2}(\tilde{D})$ such that $P_0(\hat{f}) = f$ 
and let $F \in {\cal D}_{U(1)}((\tilde{D} \times U(1))/{\bf Z}_2)$ given by 
$F([w,u]) = [\hat{f}(w), u].$  
Let $\phi \in {\cal D}_{U(1)}((\tilde{D} \times U(1))/{\bf Z}_2)$ defined by 
$\phi([w,u]) = [w,\ \exp(2 \pi v(p(w))\sqrt{-1})\, u].$ 
If $p(w) \in U, \ u \in U(1)$, we have 
%%%%%%%%%%%
\begin{eqnarray*} \lefteqn{
 (\phi^{-1} \circ F^{-1} \circ \phi \circ F)
([w,u])} \hspace{2cm} \\
&=&  (\phi^{-1} \circ F^{-1} \circ \phi)
([\hat{f}(w), u])  \\ 
&=& ( \phi^{-1} \circ F^{-1}) ([\hat{f}(w),\ 
\exp \left(2 \pi v(f(p(w)))\sqrt{-1} \, \right)\,u ])  \\
&=& \phi^{-1}([w,\ \exp \left(2\pi v(f(p(w)))\sqrt{-1}\, \right)\,u ]) \\
&=& [w,\  \exp \left(2 \pi v(p(w))\sqrt{-1} \right)^{-1}\,
\exp \left(2\pi v(f(p(w)))\sqrt{-1}\, \right)\,u ]) \\
&=& [w,\ \exp \left(2\pi \hat{a}_h(p(w))\sqrt{-1}\, \right)\,u] \\
&=& [w,\ \bar{a}_h(p(w))u] \\
&=& \gamma_*^{-1}(h)[w, u].
\end{eqnarray*}
%%%%%%%%%%%%%%%
Thus  $h =[\gamma_*^{-1}(\phi),\ \gamma_*^{-1}(F)]$ on the 
neighborhood $\pi^{-1}(U)$ of the exceptional orbit $U(1) \cdot e$. 

  Secondly assume that $\bar{a}_h(0) = -1.$  
  Since  $SL(2,{\bf R})$ is a perfect group, there exist matrices 
$A_i,B_i \in SL(2,{\bf R}) \ (1 \le i \le n)$ such that 
$$
[A_1, B_1] \cdots [A_n, B_n] = - I,
$$ 
where $I$ is the unit matrix.  
For each matrix $C \in SL(2,{\bf R})$ let $h_C \in {\cal D}_ {U(1)}(W)$ 
such that 
$$ \gamma_*^{-1}(h_C)([w,u]) = [C w, u], $$
on a neighborhood of $U(1) \cdot [0,1]$.  
Set 
$$
h' = h \circ ([h_{A_1}, h_{B_1}] \circ \cdots \circ 
[h_{A_n}, h_{B_n}])^{-1}.
$$
Since $a_{h_{- I}}(0)=-1$, we have 
 $a_{h'}(0) = a_h(0) \cdot a_{h_{- I}}(0)=1$.  
It follows from the previous argument that there exists 
$\phi \in  \left[{\cal D}_{U(1)}(W),\ {\cal D}_{U(1)}(W) \right]$ such that 
$h' = \phi$ on a neighborhood of the orbit $U(1) \cdot e$.   Hence 
$h$ also has this property.  
Thus Proposition 7.7 follows from Lemma 1 and the main theorem in \cite{A-F1}.  

\bigskip
\noindent
{\bf Proof of Theorem 7.2}\,. \quad Let $i:\ Ker \ \Theta \hookrightarrow {\cal D}_{U(1)}(W)$ be the inclusion map. 
 Then we have the following exact sequence. 
$$\begin{CD} Ker \, \Theta/\left[Ker \, \Theta,\ {\cal D}_{U(1)}(S^3) \right] @>i_*>> H_1({\cal D}_{U(1)}(S^3)) \\
  @>\Theta_*>> H_1({\cal D}(S^3/U(1)) \times U(1)/{\bf Z}_2) @>>> 1. \end{CD}$$ 
% \begin{eqnarray*}
%   Ker \, \Theta/\left[Ker \, \Theta,\ {\cal D}_ {U(1)}(S^3) \right] 
%  & \stackrel{i_*}{\to} & H_1({\cal D}_ {U(1)}(S^3)) \\ 
%  & \stackrel{\Theta_*}{\to} & H_1(D(S^3/U(1)) \times U(1)/{\bf Z}_2) \to 1.
% \end{eqnarray*}
%
By Proposition 7.7, $i_* = 0$. Then the map $\Theta_*$ is an isomorphism.  
It follows from Proposition 7.1 that 
$$H_1({\cal D}_{U(1)}(W))  \cong {\bf R} \times U(1)/{\bf Z}_2
  \cong {\bf R} \times U(1). $$
This completes the proof of Theorem 7.2.

%%%%%%%%%%%%%%%%

\vspace{5mm}

%%%%%
%\newpage
\renewcommand{\refname}

\end{document}